\documentclass{article}
\usepackage{tikz}
\usepackage{ifthen}
\usepackage{amsmath}
\usepackage{placeins}
\usepackage[colorlinks]{hyperref}
\usepackage{float}
\usepackage{caption}
\usepackage{xfp}
\usepackage{colortbl}
\providecommand{\keywords}[1]{\textbf{\textit{Keywords:}} #1}
\title{A symmetric chain decomposition of $L(5,n)$}
\author{Xiangdong Wen \\
           Wolfram Research}
\date{}
\begin{document}
\maketitle

\begin{abstract}
Stanley \cite{Stanley} conjectured that Young’s lattice $L(m, n)$ has a symmetric chain decomposition (SCD). In this paper we prove the conjecture is true for $min(m, n) = 5$ by construction a SCD for $L(5, n)$ .
\end{abstract}

\keywords{Young’s lattice, poset, symmetric chain decomposition}

\section{Introduction}
For positive integers $m$ and $n$, \textbf{Young's lattice} $L(m, n)$ refers to the partially ordered set (\textbf{poset}) of $m$-tuples $(a_1, a_2, \cdots, a_m)$, where $0 \leq a_1 \leq a_2 \leq \cdots \leq a_m \leq n$ with ordered relation
$(a_1, a_2, \cdots, a_m) \leq (b_1, b_2, \cdots, b_m )$ if $a_i \leq b_i$ for all $i=1, 2, \cdots, m$. 

The \textbf{rank} of $\overrightarrow{a}$ = $(a_1, a_2, \cdots, a_m)$ is defined as $rank(\overrightarrow{a}) = \sum_{i=1}^{m} a_i$, and a chain $\overrightarrow{v_1} < \overrightarrow{v_2} < \dots < \overrightarrow{v_k}$ in $L(m, n)$ is called \textbf{saturated} if it skips
no ranks and is called \textbf{symmetric} if $rank(\overrightarrow{v_1}) + rank(\overrightarrow{v_k}) = m n$. A symmetric chain decomposition (\textbf{SCD}) of a poset is a way of expressing it as a disjoint union of saturated symmetric chains.

Various SCDs for $L(4, n)$ and $L(3, n)$ have been constructed (\cite{West}, \cite{Lind} and \cite{Wen}).
One of the major problems in order theory is the explicit construction of SCDs for Young’s Lattice for all $m$ and $n$. In 1989, Kathy O’Hara (\cite{Kathy}, see also \cite{Doron}) constructed SCDs for the trivial extension of $L(m, n)$, in which all partitions of one rank are related to the next; but the problem remains wide open for Young’s lattice itself.

There are more than 30 types of symmetric chains in the found SCD of $L(5,n)$. To make it neat and also clear, nine types of parallel chains are presented in Section \ref{parallel chains}. We illustrate how to get symmetric chains from these parallel chains in Section \ref{symmetric chains} and finally give a computer proof in Section \ref{proof}.

\section{A list of disjointed parallel chains of $L(5, n)$} \label{parallel chains}

Here are nine types of disjointed chains of $L(5,n)$ that a SCD of $L(5,n)$ could be derived from. We name these \textbf{parallel} chains, as with varying parameters $p$,$q$ and other parameters fixed, they could form rectangles (Section \ref{symmetric chains}). 

 Parameters $i, j, k, u, p, q$ and $w$ in these tables are all non-negative integers. One line of vertical dots in these tables represents an abbreviated saturated chain with only one entry increasing. Two lines of vertical dots (in the tables $C_3, C_4, C_5$ and $C_6$) represent an abbreviated saturated chain with a zigzag path that two entries are increasing 1 alternatively. There are two helper rows with extra parameter $t$, an increasing iterator, in the zigzag path. For example, the first and the third entries are increasing by 1 alternatively at the beginning of chain $C_3$. It is also worth to mention that 
$n = 2 u + 7 + 6 j + 4 k + 3 i$ is a special case of table \ref{table:6} which does not include the first segment of the table. That is for $n = 2 u + 7 + 6 j + 4 k + 3$, parallel chains start from the darken row of table \ref{table:6}.
\\
\\
\\
\begin{table}[!tbh]
\centering
\resizebox{1 \textwidth}{!}{
\begin{tabular}{c c c c c} 
\hline
& & $  2 + 2 i + 2 j + 3 k \leq n$ \\
 \hline
 ($p$,&$ k$,&$ j + k$,&$ 1 + i + j + 2 k$,&$ 1 + i + 2 j + 2 k$)\\
&&\vdots&&\\
 ($p$,&$ k$,&$ 1 + i + j + k + p$,&$ 1 + i + j + 2 k$,&$ 1 + i + 2 j + 2 k$)\\
&\vdots&&& \\
 ($p$,&$ 1 + i + j + k$,&$ 1 + i + j + k + p$,&$ 1 + i + j + 2 k$,&$ 1 + i + 2 j + 2 k$)\\
\vdots&&&& \\
 ($1 + i + k$,&$ 1 + i + j + k$,&$ 1 + i + j + k + p$,&$ 1 + i + j + 2 k$,&$ 1 + i + 2 j + 2 k$)\\
&&&&\vdots \\
 ($1 + i + k$,&$ 1 + i + j + k$,&$ 1 + i + j + k + p$,&$ 1 + i + j + 2 k$,&$ -k + n + p$)\\
&&&\vdots& \\
 ($1 + i + k$,&$ 1 + i + j + k$,&$ 1 + i + j + k + p$,&$ -k + n$,&$ -k + n + p$)\\
&&\vdots&& \\
 ($1 + i + k$,&$ 1 + i + j + k$,&$ -j - k + n$,&$ -k + n$,&$ -k + n + p$)\\
&\vdots&&& \\
 ($1 + i + k$,&$ -1 - i - j - 2 k + n$,&$ -j - k + n$,&$ -k + n$,&$ -k + n + p$)\\
\vdots \\
 ($-1 - i - 2 j - 2 k + n$,&$ -1 - i - j - 2 k + n$,&$ -j - k + n$,&$ -k + n$,&$ -k + n + p$)\\
\hline
\end{tabular}
}
\caption{Parallel chains $C_1$, $  0 \leq p \leq k$}
\label{table:1}
\end{table}

\begin{table}[!tbh]
\centering
\resizebox{1 \textwidth}{!}{
\begin{tabular}{c c c c c} 
\hline
& & $  3 + 2 i + 2 j + 3 k \leq n$ \\
 \hline
 ($k - p$,&$ k$,&$ j + k$,&$ 1 + i + j + 2 k$,&$ 2 + i + 2 j + 2 k$)\\
&&&&\vdots \\
 ($k - p$,&$ k$,&$ j + k$,&$ 1 + i + j + 2 k$,&$ -1 - i - k + n$)\\
&&&\vdots& \\
 ($k - p$,&$ k$,&$ j + k$,&$ -1 - i - j - k + n$,&$ -1 - i - k + n$)\\
&& \vdots && \\
 ($k - p$,&$ k$,&$ -1 - i - j - k + n - p$,&$ -1 - i - j - k + n$,&$ -1 - i - k + n$)\\
& \vdots &&& \\
 ($k - p$,&$ -1 - i - j - 2 k + n$,&$ -1 - i - j - k + n - p$,&$ -1 - i - j - k + n$,&$ -1 - i - k + n$)\\
\vdots &&&& \\
 ($-1 - i - 2 j - 2 k + n$,&$ -1 - i - j - 2 k + n$,&$ -1 - i - j - k + n - p$,&$ -1 - i - j - k + n$,&$ -1 - i - k + n$)\\
&&&& \vdots \\
 ($-1 - i - 2 j - 2 k + n$,&$ -1 - i - j - 2 k + n$,&$ -1 - i - j - k + n - p$,&$ -1 - i - j - k + n$,&$ n - p$)\\
&&&\vdots & \\
 ($-1 - i - 2 j - 2 k + n$,&$ -1 - i - j - 2 k + n$,&$ -1 - i - j - k + n - p$,&$ -k + n$,&$ n - p$)\\
&& \vdots && \\
($-1 - i - 2 j - 2 k + n$,&$ -1 - i - j - 2 k + n$,&$ -1 - j - k + n$,&$ -k + n$,&$ n - p$)\\
\hline
\end{tabular}
}
\caption{Parallel chains $C_2$, $  0 \leq p \leq k$}
\label{table:2}
\end{table}

\FloatBarrier

\begin{table}[!tbh]
\centering
\resizebox{1 \textwidth}{!} {
\begin{tabular}{c c c c c} 
\hline
& & $  0 \leq u \leq 1 $ \\
& & $  2 u + 1 + 6 j + 4 k + 3 i \leq n$ \\
&&  $2 \leq t \leq i-1$ \\
 \hline
 ($2 j$,&$ i + 2 j + p$,&$ i + 2 j + 2 k + u$,&$ 2 i + 4 j + 2 k + u$,&$ 2 i + 4 j + 4 k + 2 u$)\\
 ({\color{red} $1 + 2 j$},&$ i + 2 j + p$,&$ i + 2 j + 2 k + u$,&$ 2 i + 4 j + 2 k + u$,&$ 2 i + 4 j + 4 k + 2 u$)\\
 ($1 + 2 j$,&$ i + 2 j + p$,&{\color{red} $ 1 + i + 2 j + 2 k + u$},&$ 2 i + 4 j + 2 k + u$,&$ 2 i + 4 j + 4 k + 2 u$)\\ 
{\color{red} \vdots} && {\color{red} \vdots}&&\\
 ({\color{red} $t + 2 j$},&$ i + 2 j + p$,&$  t - 1 + i + 2 j + 2 k + u$,&$ 2 i + 4 j + 2 k + u$,&$ 2 i + 4 j + 4 k + 2 u$)\\
 ($t + 2 j$,&$ i + 2 j + p$,&{\color{red} $  t + i + 2 j + 2 k + u$},&$ 2 i + 4 j + 2 k + u$,&$ 2 i + 4 j + 4 k + 2 u$)\\
{\color{red} \vdots} &&{\color{red} \vdots}&&\\
({\color{red} $i + 2 j$},&$ i + 2 j + p$,&$  2 i - 1 + 2 j + 2 k + u$,&$ 2 i + 4 j + 2 k + u$,&$ 2 i + 4 j + 4 k + 2 u$)\\
($i + 2 j$,&$ i + 2 j + p$,&{\color{red} $  2 i + 2 j + 2 k + u$},&$ 2 i + 4 j + 2 k + u$,&$ 2 i + 4 j + 4 k + 2 u$)\\
 &&&& \vdots \\
 ($i + 2 j$,&$ i + 2 j + p$,&$ 2 i + 2 j + 2 k + u$,&$ 2 i + 4 j + 2 k + u$,&$ -2 j + n$)\\
&&&\vdots& \\
 ($i + 2 j$,&$ i + 2 j + p$,&$ 2 i + 2 j + 2 k + u$,&$ -i - 2 j - 2 k + n + p - u$,&$ -2 j + n$)\\
&&\vdots&& \\
 ($i + 2 j$,&$ i + 2 j + p$,&$ -i - 2 j - 2 k + n - u$,&$ -i - 2 j - 2 k + n + p - u$,&$ -2 j + n$)\\
&\vdots&&& \\
 ($i + 2 j$,&$ -2 i - 4 j - 2 k + n - u$,&$ -i - 2 j - 2 k + n - u$,&$ -i - 2 j - 2 k + n + p - u$,&$ -2 j + n$)\\
\vdots&&& \\
($-2 i - 4 j - 4 k + n - 2 u$,&$ -2 i - 4 j - 2 k + n - u$,&$ -i - 2 j - 2 k + n - u$,&$ -i - 2 j - 2 k + n + p - u$,&$ -2 j + n$)\\
\hline
\end{tabular}
}
\caption{Parallel chains $C_3$, $0 \leq p \leq 2 k + u$}
\label{table:3}
\end{table}

\begin{table}[!tbh]
\centering
\resizebox{1 \textwidth}{!} {
\begin{tabular}{c c c c c} 
\hline
& & $  0 \leq u \leq 1 $ \\
& & $  2 u + 4 + 6 j + 4 k + 3 i \leq n$ \\
&&  $2 \leq t \leq i-1$ \\
 \hline
 ($1 + 2 j$,&$ 1 + i + 2 j + p$,&$ 1 + i + 2 j + 2 k + u$,&$ 2 + 2 i + 4 j + 2 k + u$,&$ 3 + 2 i + 4 j + 4 k + 2 u$)\\
 ($1 + 2 j$,&$ 1 + i + 2 j + p$,&{\color{red}$ 2 + i + 2 j + 2 k + u$},&$ 2 + 2 i + 4 j + 2 k + u$,&$ 3 + 2 i + 4 j + 4 k + 2 u$)\\
{\color{red} \vdots}&&{\color{red} \vdots}&& \\
 ({\color{red}$t + 2 j$},&$ 1 + i + 2 j + p$,&$ t + i + 2 j + 2 k + u$,&$ 2 + 2 i + 4 j + 2 k + u$,&$ 3 + 2 i + 4 j + 4 k + 2 u$)\\
 ($t + 2 j$,&$ 1 + i + 2 j + p$,&{\color{red}$ 1 + t + i + 2 j + 2 k + u$},&$ 2 + 2 i + 4 j + 2 k + u$,&$ 3 + 2 i + 4 j + 4 k + 2 u$)\\
{\color{red} \vdots}&&{\color{red} \vdots}&& \\
 ({\color{red}$i + 2 j$},&$ 1 + i + 2 j + p$,&$ 2 i + 2 j + 2 k + u$,&$ 2 + 2 i + 4 j + 2 k + u$,&$ 3 + 2 i + 4 j + 4 k + 2 u$)\\
 ($i + 2 j$,&$ 1 + i + 2 j + p$,&{\color{red}$ 1 + 2 i + 2 j + 2 k + u$},&$ 2 + 2 i + 4 j + 2 k + u$,&$ 3 + 2 i + 4 j + 4 k + 2 u$)\\
 ({\color{red}$1 + i + 2 j$},&$ 1 + i + 2 j + p$,&$ 1 + 2 i + 2 j + 2 k + u$,&$ 2 + 2 i + 4 j + 2 k + u$,&$ 3 + 2 i + 4 j + 4 k + 2 u$)\\
&&&&\vdots \\
 ($1 + i + 2 j$,&$ 1 + i + 2 j + p$,&$ 1 + 2 i + 2 j + 2 k + u$,&$ 2 + 2 i + 4 j + 2 k + u$,&$ -1 - 2 j + n$)\\
&&& \vdots \\
 ($1 + i + 2 j$,&$ 1 + i + 2 j + p$,&$ 1 + 2 i + 2 j + 2 k + u$,&$ -1 - i - 2 j - 2 k + n + p - u$,&$ -1 - 2 j + n$)\\
&& \vdots \\
 ($1 + i + 2 j$,&$ 1 + i + 2 j + p$,&$ -1 - i - 2 j - 2 k + n - u$,&$ -1 - i - 2 j - 2 k + n + p - u$,&$ -1 - 2 j + n$)\\
& \vdots \\
 ($1 + i + 2 j$,&$ -2 - 2 i - 4 j - 2 k + n - u$,&$ -1 - i - 2 j - 2 k + n - u$,&$ -1 - i - 2 j - 2 k + n + p - u$,&$ -1 - 2 j + n$)\\
\vdots \\
($-3 - 2 i - 4 j - 4 k + n - 2 u$,&$ -2 - 2 i - 4 j - 2 k + n - u$,&$ -1 - i - 2 j - 2 k + n - u$,&$ -1 - i - 2 j - 2 k + n + p - u$,&$ -1 - 2 j + n$)\\
 \hline
\end{tabular}
}
\caption{Parallel chains $C_4$, $0 \leq p \leq 2 k + u$}
\label{table:4}
\end{table}

\FloatBarrier

\begin{table}[!tbh]
\centering
\resizebox{1 \textwidth}{!} {
\begin{tabular}{c c c c c} 
\hline
& & $  0 \leq u \leq 1 $ \\
& & $  2 u + 4 + 6 j + 4 k + 3 i \leq n$ \\
&&  $1 \leq t \leq i - 2$ \\
 \hline
 ($2 j$,&$ 1 + i + 2 j + p$,&$ 1 + i + 2 j + 2 k + u$,&$ 2 + 2 i + 4 j + 2 k + u$,&$ 3 + 2 i + 4 j + 4 k + 2 u$)\\
 &&&& \vdots \\
 ($2 j$,&$ 1 + i + 2 j + p$,&$ 1 + i + 2 j + 2 k + u$,&$ 2 + 2 i + 4 j + 2 k + u$,&$ -1 - i - 2 j + n$)\\
 &&& \vdots \\
 ($2 j$,&$ 1 + i + 2 j + p$,&$ 1 + i + 2 j + 2 k + u$,&$ -1 - i - 2 j - 2 k + n + p - u$,&$ -1 - i - 2 j + n$)\\
 && \vdots \\
 ($2 j$,&$ 1 + i + 2 j + p$,&$ -2 - 2 i - 2 j - 2 k + n - u$,&$ -1 - i - 2 j - 2 k + n + p - u$,&$ -1 - i - 2 j + n$)\\
 & \vdots \\
 ($2 j$,&$ -2 - 2 i - 4 j - 2 k + n - u$,&$ -2 - 2 i - 2 j - 2 k + n - u$,&$ -1 - i - 2 j - 2 k + n + p - u$,&$ -1 - i - 2 j + n$)\\
 \vdots \\
($-2-2 i-4 j-4 k+n-2 u$,&$-2-2 i-4 j-2 k+n-u$,&$-2-2 i-2 j-2 k+n-u$,&$-1-i-2 j-2 k+n+p-u$,&$-1-i-2 j+n$)\\
 ($-2 - 2 i - 4 j - 4 k + n - 2 u$,&$ -2 - 2 i - 4 j - 2 k + n - u$,&{\color{red}$ -1 - 2 i - 2 j - 2 k + n - u$},&$ -1 - i - 2 j - 2 k + n + p - u$,&$ -1 - i - 2 j + n$)\\
 ($-2 - 2 i - 4 j - 4 k + n - 2 u$,&$ -2 - 2 i - 4 j - 2 k + n - u$,&$ -1 - 2 i - 2 j - 2 k + n - u$,&$ -1 - i - 2 j - 2 k + n + p - u$,&{\color{red}$ - i - 2 j + n$})\\
&&{\color{red} \vdots} &&{\color{red} \vdots} \\
 ($-2 - 2 i - 4 j - 4 k + n - 2 u$,&$ -2 - 2 i - 4 j - 2 k + n - u$,&{\color{red}$ -1 + t - 2 i - 2 j - 2 k + n - u$},&$ -1 - i - 2 j - 2 k + n + p - u$,&$ -1 + t - i - 2 j + n$)\\
 ($-2 - 2 i - 4 j - 4 k + n - 2 u$,&$ -2 - 2 i - 4 j - 2 k + n - u$,&$ -1 + t - 2 i - 2 j - 2 k + n - u$,&$ -1 - i - 2 j - 2 k + n + p - u$,&{\color{red}$        t - i - 2 j + n$})\\ 
&&{\color{red} \vdots} &&{\color{red} \vdots} \\
 ($-2 - 2 i - 4 j - 4 k + n - 2 u$,&$ -2 - 2 i - 4 j - 2 k + n - u$,&{\color{red}$ -2 - i - 2 j - 2 k + n - u$},&$ -1 - i - 2 j - 2 k + n + p - u$,&$ -2 - 2 j + n$)\\
 ($-2 - 2 i - 4 j - 4 k + n - 2 u$,&$ -2 - 2 i - 4 j - 2 k + n - u$,&$ -2 - i - 2 j - 2 k + n - u$,&$ -1 - i - 2 j - 2 k + n + p - u$,&{\color{red}$ -1 - 2 j + n$})\\
 ($-2 - 2 i - 4 j - 4 k + n - 2 u$,&$ -2 - 2 i - 4 j - 2 k + n - u$,&{\color{red}$ -1 - i - 2 j - 2 k + n - u$},&$ -1 - i - 2 j - 2 k + n + p - u$,&$ -1 - 2 j + n$)\\
 \hline
\end{tabular}
}
\caption{Parallel chains $C_5$, $0 \leq p \leq 2 k + u$}
\label{table:5}
\end{table}

\begin{table}[!tbh]
\centering
\resizebox{1 \textwidth}{!} {
\begin{tabular}{l | c c c c c} 
\hline
&& & $  0 \leq u \leq 1 $ \\
Additional conditions&& & $  2 u + 7 + 6 j + 4 k + 3 i \leq n$ \\
&&&  $1 \leq t \leq i - 2$ \\
 \hline
 $2 u + 7 + 6 j + 4 k + 3 i < n$ &($1 + 2 j$,&$ 2 + i + 2 j + p$,&$ 2 + i + 2 j + 2 k + u$,&$ 4 + 2 i + 4 j + 2 k + u$,&$ 6 + 2 i + 4 j + 4 k + 2 u$)\\
 &&&&&\vdots \\
 $2 u + 7 + 6 j + 4 k + 3 i < n$ & ($1 + 2 j$,&$ 2 + i + 2 j + p$,&$ 2 + i + 2 j + 2 k + u$,&$ 4 + 2 i + 4 j + 2 k + u$,&$ -2 - i - 2 j + n$)\\
 \rowcolor[gray]{.9}
 Chains start from this row if \\
 \rowcolor[gray]{.9}
 $2 u + 7 + 6 j + 4 k + 3 i = n$.&($1 + 2 j$,&$ 2 + i + 2 j + p$,&$ 2 + i + 2 j + 2 k + u$,&$ 5 + 2 i + 4 j + 2 k + u$,&$ -2 - i - 2 j + n$)\\
 &&&& \vdots \\
 &($1 + 2 j$,&$ 2 + i + 2 j + p$,&$ 2 + i + 2 j + 2 k + u$,&$ -2 - i - 2 j - 2 k + n + p - u$,&$ -2 - i - 2 j + n$)\\
 &&& \vdots \\
 &($1 + 2 j$,&$ 2 + i + 2 j + p$,&$ -3 - 2 i - 2 j - 2 k + n - u$,&$ -2 - i - 2 j - 2 k + n + p - u$,&$ -2 - i - 2 j + n$)\\
 && \vdots\\
 &($1 + 2 j$,&$ -4 - 2 i - 4 j - 2 k + n - u$,&$ -3 - 2 i - 2 j - 2 k + n - u$,&$ -2 - i - 2 j - 2 k + n + p - u$,&$ -2 - i - 2 j + n$)\\
 &\vdots \\
 &($-5 - 2 i - 4 j - 4 k + n - 2 u$,&$ -4 - 2 i - 4 j - 2 k + n - u$,&$ -3 - 2 i - 2 j - 2 k + n - u$,&$ -2 - i - 2 j - 2 k + n + p - u$,&$ -2 - i - 2 j + n$)\\
 &($-5 - 2 i - 4 j - 4 k + n - 2 u$,&$ -4 - 2 i - 4 j - 2 k + n - u$,&$ -3 - 2 i - 2 j - 2 k + n - u$,&$ -2 - i - 2 j - 2 k + n + p - u$,&{\color{red}$ -1 - i - 2 j + n$})\\
 &($-5 - 2 i - 4 j - 4 k + n - 2 u$,&$ -4 - 2 i - 4 j - 2 k + n - u$,&{\color{red}$ -2 - 2 i - 2 j - 2 k + n - u$},&$ -2 - i - 2 j - 2 k + n + p - u$,&$ -1 - i - 2 j + n$)\\
 &&&{\color{red} \vdots}&& {\color{red} \vdots} \\
 &($-5 - 2 i - 4 j - 4 k + n - 2 u$,&$ -4 - 2 i - 4 j - 2 k + n - u$,&$ -3 + t - 2 i - 2 j - 2 k + n - u$,&$ -2 - i - 2 j - 2 k + n + p - u$,&{\color{red}$ -1 + t - i - 2 j + n$})\\
 &($-5 - 2 i - 4 j - 4 k + n - 2 u$,&$ -4 - 2 i - 4 j - 2 k + n - u$,&{\color{red}$ -2 + t - 2 i - 2 j - 2 k + n - u$},&$ -2 - i - 2 j - 2 k + n + p - u$,&$ -1 + t - i - 2 j + n$)\\
 &&&{\color{red} \vdots}&& {\color{red} \vdots} \\
 &($-5 - 2 i - 4 j - 4 k + n - 2 u$,&$ -4 - 2 i - 4 j - 2 k + n - u$,&$ -4 - i - 2 j - 2 k + n - u$,&$ -2 - i - 2 j - 2 k + n + p - u$,&{\color{red}$ -2 - 2 j + n$})\\
 &($-5 - 2 i - 4 j - 4 k + n - 2 u$,&$ -4 - 2 i - 4 j - 2 k + n - u$,&{\color{red}$ -3 - i - 2 j - 2 k + n - u$},&$ -2 - i - 2 j - 2 k + n + p - u$,&$ -2 - 2 j + n$)\\
 &($-5 - 2 i - 4 j - 4 k + n - 2 u$,&$ -4 - 2 i - 4 j - 2 k + n - u$,&$ -3 - i - 2 j - 2 k + n - u$,&$ -2 - i - 2 j - 2 k + n + p - u$,&{\color{red}$ -1 - 2 j + n$})\\ 

\hline
\end{tabular}
}
\caption{Parallel chains $C_6$, $0 \leq p \leq 2 k + u$}
\label{table:6}
\end{table}

\FloatBarrier

\begin{table}[!tbh]
\centering
\resizebox{1 \textwidth}{!} {
\begin{tabular}{c c c c c} 
\hline
&& $u = n \mathbf{mod}{2}$, $u \in \{0,1\}$ && \\
&& $   6 + 3 u + 6 w + 6 i  + 2 k = n $ \\
 \hline
($1 + i$,&$ 2 + 2 i + p + u + 2 w$,&$ 3 + 3 i + k + u + 2 w$,&$ 4 + 4 i + k + q + 2 u + 4 w$,&$ 4 + 4 i + 2 k + 2 u + 4 w$)\\
 \vdots &&&& \\
 ($1 + i + u + 2 w$,&$ 2 + 2 i + p + u + 2 w$,&$ 3 + 3 i + k + u + 2 w$,&$ 4 + 4 i + k + q + 2 u + 4 w$,&$ 4 + 4 i + 2 k + 2 u + 4 w$)\\
 && \vdots && \\
 ($1 + i + u + 2 w$,&$ 2 + 2 i + p + u + 2 w$,&$ 3 + 3 i + k + 2 u + 4 w$,&$ 4 + 4 i + k + q + 2 u + 4 w$,&$ 4 + 4 i + 2 k + 2 u + 4 w$)\\
 ($1 + i + u + 2 w$,&$ 2 + 2 i + p + u + 2 w$,&$ 3 + 3 i + k + 2 u + 4 w$,&$ 4 + 4 i + k + q + 2 u + 4 w$,&$ 5 + 4 i + 2 k + 2 u + 4 w$)\\
\vdots \\
 ($1 + 2 i + u + 2 w$,&$ 2 + 2 i + p + u + 2 w$,&$ 3 + 3 i + k + 2 u + 4 w$,&$ 4 + 4 i + k + q + 2 u + 4 w$,&$ 5 + 4 i + 2 k + 2 u + 4 w$)\\
&&&& \vdots \\
 ($1 + 2 i + u + 2 w$,&$ 2 + 2 i + p + u + 2 w$,&$ 3 + 3 i + k + 2 u + 4 w$,&$ 4 + 4 i + k + q + 2 u + 4 w$,&$ 5 + 5 i + 2 k + 3 u + 6 w$)\\
 ($2 + 2 i + u + 2 w$,&$ 2 + 2 i + p + u + 2 w$,&$ 3 + 3 i + k + 2 u + 4 w$,&$ 4 + 4 i + k + q + 2 u + 4 w$,&$ 5 + 5 i + 2 k + 3 u + 6 w$)\\
\hline
\end{tabular}
}
\caption{Parallel chains $C_7$, $ 0 \leq p \leq k, 0 \leq q \leq p $ }
\label{table:7}
\end{table}

\begin{table}[!tbh]
\centering
\resizebox{1 \textwidth}{!} {
\begin{tabular}{c c c c c} 
\hline
&& $u = n \mathbf{mod}{2}$, $u \in \{0,1\}$   && \\
&& $   12 - 3 u + 6 w + 6 i + 2 k = n $ \\
 \hline
 ($1 + i$,&$ 4 + 2 i + p - u + 2 w$,&$ 5 + 3 i + k - u + 2 w$,&$ 8 + 4 i + k + q - 2 u + 4 w$,&$ 9 + 4 i + 2 k - 2 u + 4 w$)\\
\vdots \\
 ($3 + 2 i - u + 2 w$,&$ 4 + 2 i + p - u + 2 w$,&$ 5 + 3 i + k - u + 2 w$,&$ 8 + 4 i + k + q - 2 u + 4 w$,&$ 9 + 4 i + 2 k - 2 u + 4 w$)\\
&&&& \vdots \\
 ($3 + 2 i - u + 2 w$,&$ 4 + 2 i + p - u + 2 w$,&$ 5 + 3 i + k - u + 2 w$,&$ 8 + 4 i + k + q - 2 u + 4 w$,&$ 9 + 5 i + 2 k - 2 u + 4 w$)\\
 ($4 + 2 i - u + 2 w$,&$ 4 + 2 i + p - u + 2 w$,&$ 5 + 3 i + k - u + 2 w$,&$ 8 + 4 i + k + q - 2 u + 4 w$,&$ 9 + 5 i + 2 k - 2 u + 4 w$)\\
&& \vdots \\
 ($4 + 2 i - u + 2 w$,&$ 4 + 2 i + p - u + 2 w$,&$ 7 + 3 i + k - 2 u + 4 w$,&$ 8 + 4 i + k + q - 2 u + 4 w$,&$ 9 + 5 i + 2 k - 2 u + 4 w$)\\
&&&& \vdots \\
 ($4 + 2 i - u + 2 w$,&$ 4 + 2 i + p - u + 2 w$,&$ 7 + 3 i + k - 2 u + 4 w$,&$ 8 + 4 i + k + q - 2 u + 4 w$,&$ 10 + 5 i + 2 k - 3 u + 6 w$)\\
 \hline
\end{tabular}
}
\caption{Parallel chains $C_8$, $ 0 \leq p \leq k, 0 \leq q \leq p $}
\label{table:8}
\end{table}

\begin{table}[H]
\centering
\resizebox{1 \textwidth}{!} {
\begin{tabular}{c c c c c} 
\hline
&& $u = n \mathbf{mod}{2}$, $u \in \{0,1\}$  && \\
&& $   2 k + 3 u + 6 w = n $ \\
 \hline
 ($0$,&$ p + u + 2 w$,&$ k + u + 2 w$,&$ k + q + 2 u + 4 w$,&$ 2 k + 2 u + 4 w$)\\
\vdots \\
 ($u + 2 w$,&$ p + u + 2 w$,&$ k + u + 2 w$,&$ k + q + 2 u + 4 w$,&$ 2 k + 2 u + 4 w$)\\
 && \vdots \\
 ($u + 2 w$,&$ p + u + 2 w$,&$ k + 2 u + 4 w$,&$ k + q + 2 u + 4 w$,&$ 2 k + 2 u + 4 w$)\\
&&&&\vdots \\
 ($u + 2 w$,&$ p + u + 2 w$,&$ k + 2 u + 4 w$,&$ k + q + 2 u + 4 w$,&$ 2 k + 3 u + 6 w$)\\
 \hline
\end{tabular}
}
\caption{Parallel chains $C_9$, $0 \leq p \leq k, 0 \leq q \leq p$ }
\label{table:9}
\end{table}

\FloatBarrier

\section{Getting symmetric chains from parallel chains} \label{symmetric chains}
For fixed parameters $i, j, k$, chains of $C_1$ with parameter $p=0, 1, \cdots, k$ form a rectangle illustrated in the following figure.

\begin{figure}[H]
\centering
\resizebox{0.8 \textwidth}{!}{
\begin{tikzpicture}
\foreach \j in {1,2,3,4,5,6,7}{
  \foreach \i in {3,4,5,6,7,10,11,12,13,14}{
    \ifthenelse{ \j < 3 \OR \j > 3}{
        \fill[black] (\i,\j) circle (0.1);
        \fill[black] (\i-1,\j) circle (0.1);
        \draw[thick,->](\i-1,\j) -- (\i - 0.5,\j);
        \draw[thick,-](\i-1,\j) -- (\i ,\j);
     }{}
  }
  \foreach \i in {8}{
    \ifthenelse{ \j < 3 \OR \j > 3}{
        \draw[loosely dotted](\i-1,\j) -- (\i+1 ,\j);
     }{}
  }
}
  \draw[thick,-](2,1) -- (2,2);
  \draw[loosely dotted](2,2) -- (2,4);
  \draw[thick,-](2,4) -- (2,7);
  \draw[thick,-](14,1) -- (14,2);
  \draw[loosely dotted](14,2) -- (14,4);
  \draw[thick,-](14,4) -- (14,7);
  \draw[black] (0, 7)node[anchor=west]{$p = 0$};
  \draw[black] (0, 6)node[anchor=west]{$p = 1$};
  \draw[black] (0, 5)node[anchor=west]{$p = 2$};
  \draw[black] (0, 4)node[anchor=west]{$p = 3$};
  \draw[black] (0.5, 3)node[anchor=west]{$\vdots$};
  \draw[black] (0, 2)node[anchor=west]{$p = k-1$};
  \draw[black] (0, 1)node[anchor=west]{$p = k$};
\end{tikzpicture}
}
\caption{ A rectangle formed by parallel chains in $C_1$ for fixed $i, j, k$.}
\end{figure}
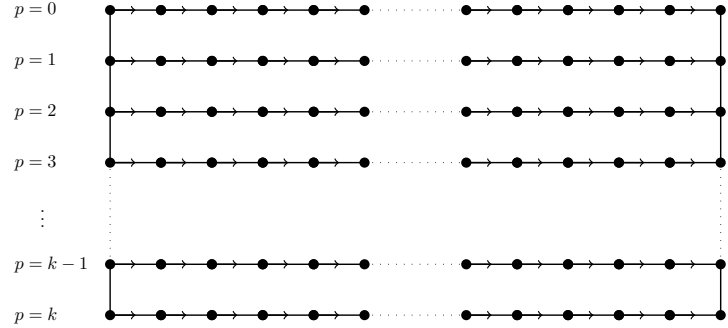

The perimeter of the rectangle contains up to 2 symmetric chains from the left up corner (the lattice point with the lowest coordinates) to the right bottom corner (the lattice point with the highest coordinates). The two corner lattice points should be in the same chain and in our cases they could be in either of the 2 border chains. All symmetric chains are obtained by taking perimeters off these rectangles recursively, illustrated in the following figure. 

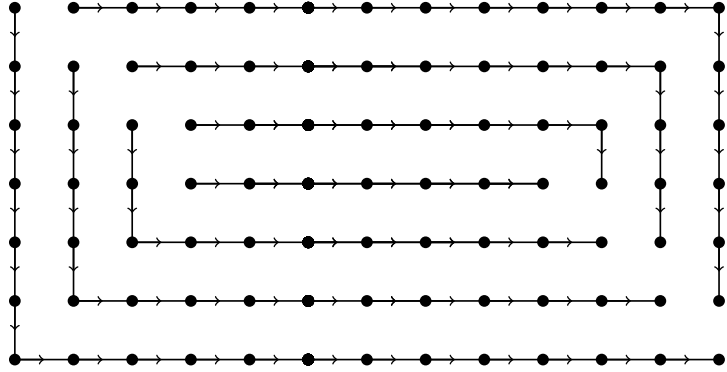
\begin{figure}[H]
\centering
\resizebox{0.8 \textwidth}{!}{
\begin{tikzpicture}
\foreach \i in {1,2,3,4,5,6,7,8,9,10,11,12,13}{
	\foreach \j in {1,2,3,4,5,6,7}{
    \fill[black] (\i,\j) circle (0.1);
    \fill[black] (6,\j) circle (0.1);
	}
	\foreach \j in {1,2,3,4,5,6,7}{
     \ifthenelse{\fpeval{\i + \j}> 8 \AND \fpeval{13 - \i} > \fpeval{7 -\j}}
	{
    \draw[thick,->](\i,\j) -- (\i + 0.5,\j);
    \draw[thick,-](\i,\j) -- (\i + 1,\j);
	}{}
     \ifthenelse{\fpeval{\i + \j}< 14 \AND \fpeval{8 - \i} < \fpeval{9 -\j}}
	{
    \draw[thick,->](\i,\j) -- (\i + 0.5,\j);
    \draw[thick,-](\i,\j) -- (\i + 1,\j);
	}{}
	\ifthenelse{\fpeval{\i - \j} > 5 \AND \fpeval{\i + \j}> 15}
	{
    \draw[thick,->](\i,\j) -- (\i ,\j-0.5);
    \draw[thick,-](\i,\j) -- (\i,\j-1);
	}{}
	\ifthenelse{\fpeval{\j - \i} >0 \AND \fpeval{\i + \j}< 9}
	{
    \draw[thick,->](\i,\j) -- (\i ,\j-0.5);
    \draw[thick,-](\i,\j) -- (\i,\j-1);
	}{}
	}
}
\end{tikzpicture}
}
\caption{Getting symmetric chains from parallel chains}
\end{figure}

We could get symmetric chains from parallel chains $C_2, C_3, \cdots, C_6$ similarly.

Parameters $(q,p)$ with conditions $0 \leq p \leq k$, $0 \leq q \leq p$ in $C_7$ ($C_8, C_9$) are actually elements in $L(2, k)$ which has a SCD
\[(t,t) \rightarrow (t,t+1) \rightarrow (t,t+2) \rightarrow  \cdots  \rightarrow  (t, k-t) \rightarrow (t+1, k-t) \rightarrow \cdots \rightarrow (k-t,k-t)\]
where $0 \leq t \leq \lfloor \frac{k}{2} \rfloor$.

Each element in a symmetric chain of $L(2, k)$ corresponds to a $L(5, n)$ chain in $C_7$ ($C_8,C_9$), and these chains together also form a rectangle illustrated in the following figure and symmetric chains could be obtained by taking perimeters off these rectangles recursively.

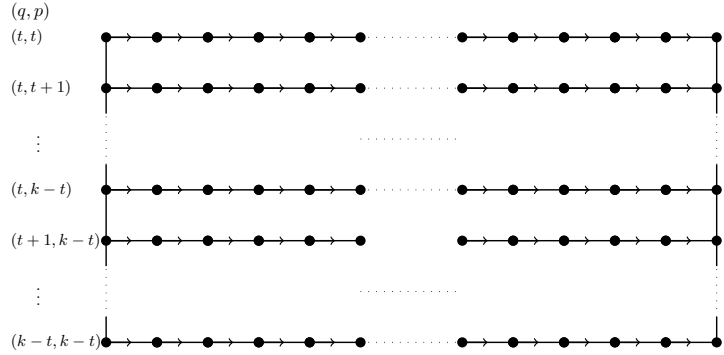
\begin{figure}[H]
\centering
\captionsetup{justification=centering}
\resizebox{0.8 \textwidth}{!}{
\begin{tikzpicture}
\foreach \j in {1,2,3,4,5,6,7}{
  \foreach \i in {3,4,5,6,7,10,11,12,13,14}{
    \ifthenelse{ \j < 2 \OR \j > 5 \OR \j = 3 \OR \j = 4}{
        \fill[black] (\i,\j) circle (0.1);
        \fill[black] (\i-1,\j) circle (0.1);
        \draw[thick,->](\i-1,\j) -- (\i - 0.5,\j);
        \draw[thick,-](\i-1,\j) -- (\i ,\j);
     }{}
   }
  \foreach \i in {8}{
    \ifthenelse{ \j < 3 \OR \j > 3}{
        \draw[loosely dotted](\i-1,\j) -- (\i+1 ,\j);
     }{}
   }
}
  \draw[thick,-](2,1) -- (2,1.5);
  \draw[loosely dotted](2,1.5) -- (2,2.5);
  \draw[thick,-](2,2.5) -- (2,4.5);
  \draw[loosely dotted](2,4.5) -- (2,5.5);
  \draw[thick,-](2,5.5) -- (2,7);
  \draw[thick,-](14,1) -- (14,1.5);
  \draw[loosely dotted](14,1.5) -- (14,2.5);
  \draw[thick,-](14,2.5) -- (14,4.5);
  \draw[loosely dotted](14,4.5) -- (14,5.5);
  \draw[thick,-](14,5.5) -- (14,7);
  \draw[black] (0, 7.5)node[anchor=west]{$(q,p)$};
  \draw[black] (0, 7)node[anchor=west]{\small $(t,t)$};
  \draw[black] (0, 6)node[anchor=west]{\small $(t,t+1)$};
  \draw[black] (0.5, 5)node[anchor=west]{\small $\vdots$};
  \draw[black] (0, 4)node[anchor=west]{\small $(t,k-t)$};
  \draw[black] (0, 3)node[anchor=west]{\small $(t+1,k-t)$};
  \draw[black] (0.5, 2)node[anchor=west]{\small $\vdots$};
  \draw[black] (0, 1)node[anchor=west]{\small $(k-t,k-t)$};
\end{tikzpicture}
}
\caption{A rectangle formed by parallel chains in $C_7$ ($C_8, C_9$) \\ for fixed $u, w, i, k$, $0 \leq t \leq \lfloor \frac{k}{2} \rfloor$}
\end{figure}

\section{Proof} \label{proof}
It is trivial to check that the chains obtained in section \ref{symmetric chains} are saturated and symmetric. 
Define $\mathbf{weight}$ of a vector $\overrightarrow{a}$ = $(a_1, a_2, \cdots , a_m)$ in $L(m, n)$ with commutative variables $x_0, x_1, x_2,\cdots, x_m$ as
\[ w(\overrightarrow{a}) = (x_0)^{n - a_m} (x_1)^{a_{m} - a_{m-1}} (x_2)^{a_{m-1} - a_{m-2}} \cdots  (x_{m-1})^{a_2 - a_1}(x_m)^{a_1} \]

For a fixed m, it is easy to see that the total weight
\[ \sum_{n=0}^{\infty} \sum_{\overrightarrow{a} \in L(m,n)} w(\overrightarrow{a}) \]
is an ordinary multivariate generating function
\[ F(x_0,x_1,\cdots, x_m) = \frac{1}{(1-x_0)(1- x_1)(1-x_2)\cdots(1-x_m)} \]

Each term in the expanded power series of $F(x_0, x_1,\cdots,x_m)$ has coefficient 1 and
corresponds to a unique vector in $L(m, n)$. On the other hand, for each vector in $L(m, n)$,
there is a unique corresponding term in the power series of $F(x_0, x_1,\cdots, x_m)$. Therefore,
to prove that each vector in $L(5, n)$ appears only once in the SCD is the same
as to prove that the total weights of the vectors in the chains is $F(x_0, x_1,x_2,x_3,x_4, x_5)$. 
As our symmetric chains are derived from parallel chains, we only need to prove

\[ \sum_{n=0}^{\infty} \sum_{i=1}^{9} \sum_{\overrightarrow{a} \in C_{i}} w(\overrightarrow{a}) = \frac{1}{(1-x_0)(1- x_1)(1-x_2)(1-x_3)(1-x_4)(1-x_5)} \]
where $C_1, C_2, \cdots, C_9$ are the parallel chains in Section \ref{parallel chains}.

This is done by computer program summing up all the weights of items in chains $C_1$, $C_2$, $\cdots$, $C_9$. The proof in Wolfram Language is available at
\href{https://www.wolframcloud.com/obj/xwen/Published/SCDL5nProof.nb}{https://www.wolframcloud.com/obj/xwen/Published/SCDL5nProof.nb}.

$SCDL5n[n]$ which gives a SCD of $L(5,n)$ written in Wolfram Language is available at \href{https://www.wolframcloud.com/obj/xwen/Published/SCDL5n.nb}{https://www.wolframcloud.com/obj/xwen/Published/SCDL5n.nb}.

\section*{Acknowledgment}
I am indebted to Dr. Kathy O’Hara for discussions which lead to finding parallel chains $C_1$ and $C_2$. 
I am indebted to Dr. Zeilberger for his help and encouragement to write up this paper.
I am also grateful to Cheng Yang for discussions and his encouragement to work on the problem.


\begin{thebibliography}{9}
\bibitem{Lind} Lindstr\"{o}m B. \emph{A partition of L(3,n) into saturated chains}, European J. Comb,
61-63,1(1980).

\bibitem{Kathy} O’Hara K.M. \emph{Unimodality of Gaussian coefficients: a constructive proof},
J. Combin. Theory Ser. A 29-52, 53 (1990).

\bibitem{Stanley} Stanley, R. \emph{Weyl groups, the hard Lefschetz theorem, and the Sperner property}, 
SIAM J. Algebr. Discrete Methods 1, 168–184, 1980.

\bibitem{Wen} Wen X., \emph{Computer-generated symmetric chain decompositions for L(4, n)
and L(3, n)} Adv. Appl. Math., v. 33 (2004) 409-412.

\bibitem{West} West D.B. \emph{A symmetric chain decomposition of L(4,n)}, European J. Comb,
379-383,1(1980).

\bibitem{Doron} Zeilberger D. \emph{Kathy O’Hara’s Constructive proof of the Unimodality of the
Gaussian Polynomials}, American Mathematical Monthly, 590-602,96(1989)

\end{thebibliography}
\end{document}